
\baselineskip=14pt
\parskip=10pt
\def\halmos{\hbox{\vrule height0.15cm width0.01cm\vbox{\hrule height
  0.01cm width0.2cm \vskip0.15cm \hrule height 0.01cm width0.2cm}\vrule
  height0.15cm width 0.01cm}}

\magnification=\magstephalf

\def\1{{\overline{1}}}
\def\2{{\overline{2}}}
\parindent=0pt
\overfullrule=0in

\def\frac#1#2{{#1 \over #2}}

\bf
\centerline
{
Computerizing  the  Andrews-Fraenkel-Sellers Proofs on the Number 
}
\centerline
{
of m-ary partitions mod m (and doing MUCH more!)
}

\rm
\bigskip
\centerline
{\it By Shalosh B. EKHAD and Doron ZEILBERGER}
\bigskip

{\bf VERY IMPORTANT} 

As in all our joint papers, the main point is not the article, but the accompanying Maple package, 
that for the present article happens to be {\tt AFS.txt}. It
may be downloaded, free of charge, from the webpage of this article

{\tt http://www.math.rutgers.edu/\~{}zeilberg/mamarim/mamarimhtml/afs.html} \quad ,

where the readers can also find sample input and output files, that they are welcome to extend using their
own computers.

{\bf A Concise Rendition of the Work of Andrews, Fraenkel, and Sellers on the number of m-ary partitions mod m}

In two delightful articles ([AFS1],[AFS2]), George Andrews, Aviezri Fraenkel, and Jim Sellers prove two results
that are trivially equivalent to the following two propositions.

{\bf Proposition A} ([AFS1], Theorem 2.1) Fix a positive integer $m$ larger than $1$.
Let 
$$
B(q)= \sum_{i=0}^{\infty} \, b(i) \, q^i \quad,
$$
be the {\it unique} formal power series satisfying the {\it functional equation}
$$
B(q)=\frac{1}{1-q} \, B(q^m) \quad,
\eqno(FE1)
$$
and let
$$
b_0(n):= \, b(m\,n) \quad .
$$
Then $b_0(n) ({\bmod \,\, m})$ can be computed in {\it logarithmic time} via the recurrence
(recall that, thanks to Euclid, every integer $n$ can be written as $n=mi+j$, where
$i$ is the {\it quotient} and $j$ ($0 \leq j <m$) is the {\it remainder})
$$
b_0(mi+j) \equiv (j+1) \, b_0(i) \,\,\, ({\bmod \,\, m}) \quad .
$$

{\bf Proposition B} ([AFS2], Theorem 2.1) Fix a positive integer $m$ larger than $1$.
Let 
$$
C(q)= \sum_{i=0}^{\infty} \, c(i)\, q^i \quad,
$$
be the {\it unique} formal power series satisfying the {\it functional equation}
$$
C(q)= 1 \, + \, \frac{q}{1-q} \, C(q^m) \quad,
\eqno(FE2)
$$
and let
$$
c_1(n):= \, c(mn \, + \, 1) \quad .
$$
Then $c_1(n) \,\, ({\bmod \,\, m})$ can be computed in {\it logarithmic time} via the recurrence
(recall, that thanks to Euclid, every integer $n$ can be written as $n=mi+j$, where
$i$ is the {\it quotient} and $j$ ($0 \leq j <m$) is the {\it remainder})
$$
c_1(mi+j) \, \equiv  \, 1+ j \, c_1 (i) \, \, ({\bmod \,\, m}) \quad .
$$

{\bf Comment}: The original statement in [AFS2] regarded $c_0(n)=c(mn)$ (rather than $c_1(n)=c(mn+1)$) and 
was rather complicated. Since it is
readily seen that $c_0(n)=c_1(n-1)$, our formulation is  simpler.

The basic {\bf idea} behind the proofs of these two propositions is {\it brilliant}, but the way the proofs are presented there are
unnecessarily long. We will first present new renditions of their nice proofs, in a form that would make them
amenable to formulate a general {\bf algorithm} that can handle many other cases.

{\bf Definition}: The $m$-sections of a formal power series, $f(q)$,  are the unique formal power series
$f_0(q), \dots, f_{m-1}(q)$ such that
$$
f(q)=\sum_{i=0}^{m-1} \, q^i \, f_i (q^m) \quad .
$$

{\bf Remark}: It is an elementary exercise in {\it Linear Algebra} to prove that if $f(q)$ is a rational function, so
are all the $f_i(q)$'s and it is routine to find them. Of course, one can use ``averaging'' over roots of unity,
but it is not necessary.

{\bf Proof of Proposition A} ([AFS1], streamlined by DZ) : Let us $m$-sect both sides of the defining functional
equation and extract the part consisting of powers that are multiples of $m$
$$
B_0(q^m)=\left [ \left ( \frac{1}{1-q} \right )_0 (q^m) \right ] \, \cdot \, B(q^m) \quad .
\eqno(FE1')
$$
But
$$
\left ( \frac{1}{1-q} \right)_0 (q^m)= \sum_{i=0}^{\infty} q^{im} = \frac{1}{1-q^m} \quad .
$$
Hence
$$
B_0(q^m)=\frac{1}{1-q^m} \, B(q^m) \quad .
\eqno(FE1')
$$
Replacing $q^m$ by $q$ we get
$$
B_0(q)=\frac{1}{1-q} \, B(q) \quad .
\eqno(FE1'')
$$
Plugging $B(q)=(1-q)B_0(q)$ into $(FE1)$, we get
$$
(1-q)B_0(q) = \frac{1}{1-q} \, (1-q^m)B_0(q^m) \quad ,
$$
and we get a {\bf functional equation} for $B_0(q)$
$$
B_0(q)= \frac{1-q^m}{(1-q)^2} \, B_0(q^m) \quad .
$$

Now take the {\bf partial fraction decomposition} 
$$
\frac{1-q^m}{(1-q)^2} = \frac{m}{1-q} \, + \, \sum_{j=0}^{m-1} (j+1-m) q^j \quad .
$$
Hence
$$
B_0(q) \, \equiv  \left ( \sum_{j=0}^{m-1} (j+1) q^j \right )\, B_0(q^m) \quad ({\bmod \,\,m}) \quad,
$$
that is equivalent to the stated recurrence, by extracting the coefficient of $q^{mi+j}= q^j \cdot (q^m)^i $ from both sides. \halmos

{\bf Proof of Proposition B} ([AFS2], streamlined by DZ) : Now we take the $f_1(q)$ part of both sides of $(FE2)$, getting
$$
q\,C_1(q^m)= (1)_1+ q \left [ \, \left ( \frac{q}{1-q} \right)_1(q^m) \, \right ]\,\cdot \, C (q^m) \quad .
\eqno(FE2')
$$
Since $(1)_1=0$ and  $\left ( \frac{q}{1-q} \right)_1=\sum_{i=0}^{\infty} q^i \, = \, \frac{1}{1-q}$, we get
$$
C_1(q^m)= \frac{1}{1-q^m} \,  C(q^m) \quad .
$$
Replacing $q^m$ by $q$, we get
$$
C_1(q)= \frac{1}{1-q} \,  C(q) \quad ,
$$
hence
$$
C(q)= (1-q) \,  C_1 (q) \quad ,
$$
that leads to the functional equation for $C_1(q)$:
$$
C_1(q)= \frac{1}{1-q} \, + \, \frac{q(1-q^m)}{(1-q)^2} \, C_1(q^m) \quad .
$$
Now take the {\bf partial fraction decomposition} 
$$
\frac{q(1-q^m)}{(1-q)^2} = \frac{m}{1-q} \, + \, \sum_{j=0}^{m-1} (j-m) q^j \quad .
$$
Hence, taking it modulo $m$, we get,
$$
C_1(q) \, \equiv  \frac{1}{1-q} \, + \, \left ( \sum_{j=0}^{m-1} j\, q^j \right )\, C_1(q^m) \quad ({\bmod \,\,m}) \quad,
$$
that is equivalent to the stated recurrence, by extracting the coefficient of $q^{mi+j}= q^j \cdot (q^m)^i $ from both sides. \halmos

{\bf Let's Generalize!}

The Andrews-Fraenkel-Sellers method of proof (once compactified, as above) suggests an {\bf algorithm} that
tries to find poly-logarithmic-time (i.e. polynomial in the bit-size) schemes for the congruence class modulo $m$ for the $i$-part in the $m$-section
of {\bf any} formal power series satisfying a functional equation of the form
$$
f(q)=S(q) + R(q) f(q^m) \quad ,
$$
for {\bf any} given rational functions $S(q)$ and $R(q)$ 
(whose denominators do not vanish at $q=0$ so they are bona-fide formal power series, and $R(0)S(0)=0$,
in order there to be a solution)
and {\bf any} $i$ between $0$ and $m-1$. Of course, now $m$, and $i$ have
to be {\it specific}, i.e. {\it numeric}, but if in luck, one can easily detect the general pattern in $m$
(that can be conjectured and {\it proved} automatically).

The reason that things worked out so well in Propositions A and B above was that after the
partial fraction decomposition, taking it modulo $m$, the rational function part {\bf disappeared}
and we were left with a {\bf polynomial} on the right side of the functional equation for
$B_0(q)$ and $C_1(q)$ modulo $m$. But why not try, and look for more miracles?

{\bf Algorithm} 

{\bf Inputs}

$\bullet$ A positive integer $m>1$ and a non-negative integer $i$, $0\leq i <m$ \quad .

$\bullet$ Rational functions $S(q)$ and $R(q)$ whose denominators do not vanish at $q=0$, and $R(0)S(0)=0$.

Let $F(q)$ be the unique formal power series that satisfies
$$
F(q)=S(q) \, + \, R(q) \, F(q^m) \quad.
\eqno(FE3)
$$

Let $F_i(q)$ be the $i$-th part in the $m$-section
$$
F(q)=\sum_{i=0}^{m-1} q^i F_i(q^m) \quad .
$$

{\bf Output}: A rational function $E(q)$ and polynomial $P(q)$ such that $F_i(q)$ satisfies the fast functional equation, modulo $m$
$$
F_i(q) \, \equiv  E(q) \,+\, P(q) \, F_i (q^m) \quad  ( {\bmod \,\,m} ).
$$
or else FAIL.

[Note that one can compute the coefficients of rational functions, modulo $m$, in logarithmic time.]

{\bf Description}

 First $m$-sect the functional equation $(FE3)$ and extract the $i$-part
$$
q^i \, F_i(q^m)= q^i S_i(q^m) \, + \, q^i R_i(q^m) \, F(q^m) \quad.
\eqno(FE3')
$$
Note that the computer can easily compute the rational functions $S_i(q),R_i(q)$, by the
$m$-section procedure (implemented in procedure {\tt mSectR} in {\tt BFF.txt}).

Divide by $q^i$ and replace $q^m$ by $q$, getting
$$
 F_i(q)= S_i(q) \, + \,  R_i(q) \, F(q) \quad.
$$
Hence
$$
F(q)= \frac{F_i(q)-S_i(q)}{R_i(q)} \quad .
$$
Plug this into $(FE3)$, and get a brand-new functional equation for $F_i(q)$. Suppose it is
$$
F_i(q)=A(q) \, + \, G(q) \, F_i(q^m) \quad,
\eqno(FE4)
$$
for some rational functions $A(q)$ and $G(q)$ that the computer can easily find.
To wit:
$$
A(q)= 
\frac{
S \left( q \right)  
R_i \left( q \right) 
 R_i \left( {q}^{m} \right) 
-R \left( q \right) 
 R_i \left( q \right)  
S_i \left( {q}^{m} \right) 
+ R_i \left( {q}^{m} \right) 
S_i \left( q \right) 
 R_i \left( {q}^{m} \right)} 
{ R_i \left( {q}^{m} \right) }
\quad ,
$$
$$
G(q)=
\frac{R \left( q \right)  R_i \left( q \right) }
{ R_i \left( {q}^{m} \right) } \quad .
$$

Next perform the {\bf partial fraction decomposition} of $G(q)$:
$$
G(q)= PurelyRationalPart(q) + PolynomialPart(q) \quad .
$$
If a {\bf miracle} {\it does} happen, in other words, $PurelyRationalPart(q)$ is a multiple of $m$, then, we get
$$
F_i(q) \equiv A(q) \, + \, PolynomialPart(q) \, F_i(q^m) \quad ( {\bmod \,\,m} ) \quad,
$$
and we get a poly-logarithmic time way to compute the coefficients of $F_i(q)$ modulo $m$.
If the miracle does not happen, then we return {\tt FAIL}.

{\bf Another Miracle}

{\bf Proposition C}

Let $c(n)$ be the number of partitions of $n$ into parts that are either powers of $m$ or twice powers of $m$, so
the generating function
$$
C(q)= \sum_{n=0}^{\infty} c(n)q^n = \prod_{i=0}^{\infty} \frac{1}{(1-q^{m^i})(1-q^{2m^i})} \quad,
$$
satisfies the functional equation
$$
C(q) = \frac{1}{(1-q)(1-q^2)} C(q^m) \quad .
$$
Let $d(n):=\,c(mn+m-1)=c(m(n+1)-1)$, and let $D(q)=\sum_{n=0}^{\infty} d(n)q^n$ be its generating function.

Then  $d(n) \, (\,{\bmod \,\,m}\,)$ can be computed in logarithmic time, via
$$
D(q) \equiv  Nes (q) D(q^m) \quad ({\bmod \,\, m}) \quad,
$$
where, $Nes(q)$ is the polynomial of $q$ defined as follows. If $m$ is odd, then
$$
Nes(q) \, := \,(1+q) \sum_{j=0}^{m-2} {{j+2} \choose {2}} q^{2j} \quad ({\bmod \,\, m}) \quad ,
$$
while if $m$ is divisible by $4$, then
$$
Nes(q) \, := \, \sum_{j=0}^{m-3} A002623(j) (q^j+ q^{2m-4-j}) ({\bmod \,\, m}) \quad,
$$
where (see [S]) $A002623(j)=\lfloor \frac{(j+2)(j+4)(2j+3)}{24} \rfloor$.

There are no miracles when $m \equiv 2 \,\, ({\bmod \,\, 4})$.

{\bf Infinitely  more miracles}

The reader can find a few more miracles in the sample output files given in the
front of this article, mentioned above:

{\tt http://www.math.rutgers.edu/\~{}zeilberg/mamarim/mamarimhtml/afs.html}  \quad ,

and potentially infinitely more, by playing with the accompanying Maple package 
{\tt AFS.txt} available from the above page, or directly from

{\tt http://www.math.rutgers.edu/\~{}zeilberg/tokhniot/AFS.txt}  \quad .

{\bf Conclusion}

Thomas Edison said that genius is $\%1$ {\it inspiration} and $\%99$ {\it perspiration}.
Now that we have computers, they can do the perspiration part for us, but we need
{\it meta-inspiration}, {\it meta-geniuses},
and {\it meta-perspiration}, to teach the human inspiration  to our
silicon colleagues. Sooner or later, computers will also do the inspiration part, but
let humans enjoy the remaining fifty (or whatever) years left for them, and focus on
inspiration, meta-inspiration, and meta-perspiration,
and leave the {\it actual} perspiration part to their much faster-
and {\bf much more reliable}- machine friends.

{\bf Acknowledgment}

We wish to thank Aviezri Fraenkel and George Andrews for stimulating conversations.

{\bf References}

[AFS1] G.E. Andrews, A. S. Fraenkel, and J. A. Sellers, {\it Characterizing the number if $m$-ary partitions modulo $m$},
American Mathematical Monthly (2016), forthcoming.

[AFS2] G.E. Andrews, A. S. Fraenkel, and J. A. Sellers, {\it $m$-ary partitions with no gaps: A characterization modulo $m$},
Discrete Mathematics {\bf 339}(2016), 283-287.

[S] N.J.A. Sloane, {\it The On-Line Encyclopedia of Integer Sequences}, {\tt www.oeis.org} .

\bigskip
\bigskip
\hrule
\bigskip
Doron Zeilberger, Department of Mathematics, Rutgers University (New Brunswick), Hill Center-Busch Campus, 110 Frelinghuysen
Rd., Piscataway, NJ 08854-8019, USA. \hfill \break
zeilberg at math dot rutgers dot edu \quad ;  \quad {\tt http://www.math.rutgers.edu/\~{}zeilberg/} \quad .
\bigskip
\hrule
\bigskip
Shalosh B. Ekhad, c/o D. Zeilberger, Department of Mathematics, Rutgers University (New Brunswick), Hill Center-Busch Campus, 110 Frelinghuysen
Rd., Piscataway, NJ 08854-8019, USA.
\bigskip
\hrule

\bigskip
Exclusively published in The Personal Journal of Shalosh B. Ekhad and Doron Zeilberger  \hfill \break
({ \tt http://www.math.rutgers.edu/\~{}zeilberg/pj.html})
and {\tt arxiv.org} \quad . 
\bigskip
\hrule
\bigskip
{\bf Nov. 20, 2015}

\end